\documentclass[english]{article}
\usepackage{amssymb}
\usepackage{amsmath}
\usepackage{babel}

\def\thebibliography#1{\section*{References}\list
  {[\arabic{enumi}]}{\settowidth\labelwidth{[#1]}\leftmargin\labelwidth
    \advance\leftmargin\labelsep
    \usecounter{enumi}}
    \def\newblock{\hskip .11em plus .33em minus -.07em}
    \sloppy
    \sfcode`\.=1000\relax}
\newcommand{\refbook}[3]{{\sc #1}{\em\ #2}{\ #3}}
\newcommand{\refer}[5]{{\sc #1}{\ #2}{\em\ #3}{\bf\ #4}{\ #5}}

\newtheorem{lem}{Lemma}[section]
\newtheorem{cor}[lem]{Corollary}
\newtheorem{teo}[lem]{Theorem}
\newtheorem{os}[lem]{Remark}

\newtheorem{prop}[lem]{Proposition}

\newcommand{\qed}{\thinspace\null\nobreak\hfill\hbox{\vbox{\kern-.2pt\hrule
 height.2pt depth.2pt\kern-.2pt\kern-.2pt \hbox to2.5mm{\kern-.2pt\vrule
 width.4pt \kern-.2pt\raise2.5mm\vbox to.2pt{}\lower0pt\vtop
 to.2pt{}\hfil\kern-.2pt \vrule
 width.4pt \kern-.2pt}\kern-.2pt\kern-.2pt\hrule height.2pt depth.2pt
 \kern-.2pt}}\par\medbreak}

\newcommand{\R}{\mathbb{R}}

\newcommand{\C}{\mathbb{C}}

\newcommand{\N}{\mathbb{N}}

\newcommand{\eps}{\varepsilon}

\newcommand{\ov}{\overline}

\newcommand{\ds}{\displaystyle}
\newcommand{\supp}{\emph{supp\,}}

\date{}
\begin{document}

\title{Elliptic operators with unbounded diffusion coefficients in $L^p$ spaces}
\author{G. Metafune, C. Spina\thanks{Dipartimento di Matematica ``Ennio De
Giorgi'', Universit\`a del Salento, C.P.193, 73100, Lecce, Italy.
e-mail: giorgio.metafune@unisalento.it, chiara.spina@unisalento.it}}
\maketitle
\begin{abstract}
In this paper we prove that,  under suitable assumptions on $\alpha>0$, the operator $L= (1+|x|^\alpha)\Delta$ admits realizations generating contraction or analytic semigroups in $L^p(\R^N)$. For some values of $\alpha$, we also explicitly characterize the domain of $L$. Finally, some informations about the location and composition of the spectrum are given.

\end{abstract}

\par


\section{Introduction}
In this paper we focus our attention on a class of elliptic
operators with unbounded diffusion coefficients. We deal with
operators of the form
\begin{equation} \label{defL}
Lu=(1+|x|^\alpha)\Delta u,
\end{equation}
for positive values of $\alpha$,  on $L^p=L^p(\R^N, dx)$ with
respect to the Lebesgue measure. The case $\alpha\leq 2$ has been
already investigated in literature and for this reason we shall
assume $\alpha >2$ throughout the paper, even when some argument
easily extends to lower values of $\alpha$. We refer to
\cite{for-lor} where it is proved that the operator above generates
a strongly continuous and analytic semigroup in $L^p$ and in spaces
of continuous functions. For $1<p<\infty$ an explicit description
follows from the a-priori estimates
$$\|(1+|x|^\alpha)D^2u\|_p\leq C(\|u\|_p+\|(1+|x|^\alpha)\Delta u\|_p.$$
Similar estimates hold for a mor general class of operators, they
can be deduced by some weigthed norm inequalities for
Cald\'{e}ron-Zygmund singular integrals. Muckenhoupt and Wheeden for
example (see \cite{muck-wheed} or \cite{stein-weiss}) proved that
estimates of the form
$$\|aD^2 u\|_p\leq C\|a\Delta u\|_p$$
are true  for a weight $a$ in some suitable Muckenhoupt classes. In
particular the estimates above imply that
\begin{equation}   \label{est}
\||x|^\alpha D^2 u\|_p\leq C\||x|^\alpha\Delta u\|_p
\end{equation}
and
$$\|(1+|x|^\alpha)D^2u\|_p\leq C(\|u\|_p+\|(1+|x|^\alpha)\Delta u\|_p$$
for $0<\alpha<\frac{N}{p'}$ where $p'$ is the conjugate exponent of $p$.\\
Similar estimates follow also by \cite{kree} where the author proved
that certain singular integrals
are convolution operators in weighted $L^p$ spaces for the weight $1+|x|^\alpha$, $-\frac{N}{p'}<\alpha<\frac{N}{p'}$. \\
We will prove that for certain values of $\alpha>2$  the operator
above admits realizations generating analytic semigroups in $L^p$
for $1<p<\infty$.
Moreover for some values of $\alpha$ we will give also an explicit description of the domain by proving some a-priori estimates.\\
The starting point is a generation result of strongly continuous
semigroups in spaces of continuous functions. It is known (see
\cite[Example 7.3]{met-wack}) that, if $N\geq 3$, $\alpha>2$, the
operator generates a strongly continuous semigroup in $C_0(\R^N)$,
it has also been proved that both the semigroup and the resolvent
are compact.

{\bf Notation.} We use $L^p$ for $L^p(\R^N, dx)$, where this latter
is understood with respect to the Lebesgue measure. $C_b(\R^N)$ is
the Banach space of all continuous and bounded functions in $\R^N$,
endowed with the sup-norm, and $C_0(\R^N)$ its subspace consisting
of all continuous functions vanishing at infinity. $C_c^\infty(\R^N)$
denotes the set of all $C^\infty$ functions with compact support.

\section{Solvability in spaces of continuos functions}
The solvability of elliptic and parabolic problems associated to $L$
in $L^p$ depends on $\alpha, p, N$. However, these restricions are
not necessary in $C_b(\R^N)$ for a larger class of operators.
Following \cite{met-wack}, we recall the main results in spaces of
continuos functions which will be useful for comparison throughout
the paper

Let $A$ be a second
order elliptic partial differential operator of the form
$$Au(x)=\sum_{i,j=1}^Na_{ij}(x)D_{ij}u(x)+\sum_{i=1}^N F_i(x)D_iu(x)\quad\quad\quad x\in\R^N$$
under the following hypotheses on the coefficients: $a_{ij}=a_{ji}$,
$a_{ij},\ F_i$ are real-valued   locally H\"older continuous
functions of exponent $0<\alpha<1$  and  the matrix $(a_{ij})$
satisfies the ellipticity condition
$$\sum_{i,j=1}^Na_{ij}(x)\xi_i\xi_j\geq \lambda(x)|\xi|^2$$
for every $x,\ \xi\in\R^N$, with $\inf_K\lambda(x)>0$ for every compact $K\subset\R^N$. The operator $A$ is locally uniformly elliptic, that is uniformly elliptic on every compact subset of $\R^N$.\\
We endow $A$ with its maximal domain in $C_b(\R^N)$ given by
$$D_{max}(A)=\{u\in C_b(\R^N)\cap W^{2,p}_{loc}(\R^N)\quad\textrm{for all}\quad p<\infty:\ Au\in C_b(\R^N)\}.$$
The main interest is in the existence of (spatial) bounded
solutions of the parabolic problem
\begin{equation}     \label{problem}
\left\{
\begin{array}{ll}
u_t(t,x)=Au(t,x)& x\in\R^N,\ t>0, \\
u(0,x)=f(x)   &x\in\R^N \end{array}\right.
\end{equation}
with initial datum $f\in C_b(\R^N)$. The unbounded interval
$[0,\infty[$ can be changed to  any bounded
$[0,T]$ without affecting the results. Since the coefficients can be
unbounded, the classical theory does not apply and existence and
uniqueness for (\ref{problem}) are not clear. Quite surprisingly,
existence is never a problem as stated in the following theorem.
\begin{teo}   \label{exist-sem}
There exists a positive  semigroup $(T(t))_{t\geq 0}$ defined in
$C_b(\R^N)$ such that, for any $f\in C_b(\R^N)$, $u(t,x)=T(t)f(x)$
belongs to the space
$C_{loc}^{1+\frac{\alpha}{2},2+\alpha}((0,+\infty)\times \R^N)$, is
a bounded solution of the following differential equation
$$u_t(t,x)=\sum_{i,j=1}^Na_{ij}(x)D_{ij}u(t,x)+\sum_{i=1}^N F_i(x)D_iu(t,x)$$
and satisfies
$$\lim_{t \to 0}u(t,x)=f(x)$$
pointwise.
\end{teo}
When $f \in C_0(\R^N)$, then $u(t,\cdot) \to f$ uniformly as $t \to
0$. This, however, does not mean that $T(t)$ is strongly continuous
on $C_0(\R^N)$ since this latter need not to be preserved by the
semigroup

The idea of the proof is to take an increasing sequence of balls
filling the whole space and, in each of them, to find a solution  of
the parabolic problem associated with the operator. Then the
sequence of solutions so obtained is proved to converge to a
solution of the problem in $\R^N$. More precisely,  let us fix a
ball $B_\rho=B_{\rho}(0)$ in $\R^N$ and consider the problem
 \begin{equation}     \label{problem-bound}
\left\{
\begin{array}{ll}
u_t(t,x)=Au(t,x) & x\in B_\rho,\,\ t>0, \\
u(t,x)=0    &x\in \partial B_\rho,\ t>0\\
u(0,x)=f(x)   &x\in\R^N. \end{array}\right.
\end{equation}
Since the operator $A$ is uniformly elliptic and the coefficients
are bounded in $B_\rho$, there exists a unique solution $u_\rho$ of
problem (\ref{problem-bound}).  The next step consists in letting
$\rho$ to infinity in order to define the semigroup as\-so\-cia\-ted
with $A$ in $\R^N$. By using the parabolic maximum principle, it is
possible to prove that the sequence $u_\rho$ increases with $\rho$
when $f \ge 0$ and is uniformly bounded by the sup-norm of $f$. In
virtue of this monotonicity and since a general $f$ can be written as
$f=f^+-f^-$, the limit
$$T(t)f(x):=\lim_{\rho\to\infty} u_\rho(t,x)$$
is well defined for  $f\in C_b(\R^N)$ and one shows all relevant
properties, using the interior Schauder estimates.

It is worth-mentioning that also the resolvent of $A$, namely
$(\lambda-A)^{-1}$, is, for positive $\lambda$, the limit as $\rho
\to \infty$ of the corresponding resolvents in the balls $B_\rho$.
The construction then shows that, for positive $f \in C_b(\R^N)$ and
$\lambda >0$, both the semigroup $T(t)f$ (and the resolvent
$(\lambda-A)^{-1}f$) select in a linear way the minimal solution
among all bounded solutions of (\ref{problem}) (of $\lambda u-Au=f)$.
For this reason, from now on, the semigroup $T(t)$ will be called
the minimal semigroup associated to $A$ and will be denoted by
$T_{\min}(t)$. Its generator $(A,D)$, where $D \subset D_{\max}(A)$,
will be denoted by $A_{\min}$

In contrast with the existence, the uniqueness is not guaranteed, in
general, and relies on the existence of suitable Lyapunov functions.
We do not deal here with such a topic and refer again to
\cite{met-wack}. We only point out that uniqueness holds if and only
if $D=D_{\max }(A)$, i.e. when $A_{\min}$ coincides with
$(A,D_{\max}(A))$.

Let us specialize to our operator $L$.

\begin{prop}     \label{proprietaL}
Let $L= (1+|x|^\alpha)\Delta$.
\begin{itemize}
\item[(i)]
If $\alpha\leq 2$, the semigroup preserves $C_0(\R^N)$ and neither
the semigroup nor the resolvent are compact.\\
\item[(ii)] If $\alpha>2$ and $N=1,\ 2$, the semigroup is generated by $(A, D_{max}(A))$,
$C_0(\R^N)$ and $L^p$ are not preserved by the semigroup and the
resolvent and
both the semigroup and the resolvent are compact.\\
\item[(iii)] If $\alpha >2$,  $N\geq 3$, then the semigroup is generated by $(A,
D_{max}(A))\cap C_0(\R^N)$, the resolvent and the semigroup map
$C_b(\R^N)$ into $C_0(\R^N)$ and are compact.
\end{itemize}
\end{prop}
See (\cite[Example 7.3]{met-wack}). In particular (ii) will imply
that if $\alpha >2$ and $N=1,2$, problem (\ref{problem}) cannot be solved
in $L^p$. Observe also that (iii) and the discussion above show that
$(T_{min}(t))_{t \ge 0}$ is strongly continuous on $C_0(\R^N)$.

\section{Preliminary considerations in $L^p$}
We consider the operator $\hat{L}_p=(L,\hat{D}_p)$ on any domain
$\hat{D}_p$ contained in the  maximal domain in $L^p(\R^N)$ defined
by
\begin{equation} \label{massimale}
D_{p,max}(L)=\{u\in L^p\cap W^{2,p}_{loc}:\ Lu\in L^p\}.
\end{equation}
Note that $D_{p,max}(L)$ is the analogous of $D_{max}(L)$ for
$p<\infty$.
 We are interested in solvability of elliptic and parabolic
problems associated to $L$. We show that for certain values of $p$
the equation
$$
\lambda u-Lu=f
$$
is not solvable  in $L^p(\R^N)$ for positive $\lambda$. \\

In the following proposition we show that functions in $D_{p,max}(L)$
are globally in $W^{2,p}$.

\begin{prop} \label{w2p}
$$D_{p,max}(L)=\{u \in W^{2,p}: (1+|x|^\alpha)\Delta u \in L^p\}.$$
\end{prop}
{\sc Proof.} It is clear that the right hand side is included in the
left one. Conversely, if $u \in D_{p,max}(L)$, then $u,\ \Delta u \in
L^p$ and we have to show that $u \in W^{2,p}$. Let $v \in W^{2,p}$
be such that $v-\Delta v=u-\Delta u$. Then $w=u-v \in L^p$ solves
$w-\Delta w=0$. Since $w$ is a tempered distribution, by taking the
Fourier transform it easily follows that $w=0$, hence $u=v$. \qed

 The next lemma shows that the resolvent operator in
$L^p$, if it exists, is a positive operator.
\begin{lem}        \label{positiv}
Suppose that $\lambda\in\rho(\hat{L}_p)$ for some $\lambda\geq 0$.
Then for every $0\le f\in L^p$,
$$(\lambda-\hat{L}_p)^{-1}f\geq 0.$$
\end{lem}
{\sc Proof.} By density we may assume that $0\leq f\in
C^\infty_c(\R^N)$. Suppose first that $\lambda >0$. We set
$u=(\lambda-L_p)^{-1}f$. Suppose $\supp{f}\subset B(R)$. Then $u$
satisfies
$$\lambda u-Lu=f$$ in $B(R)$ and $$\lambda u-Lu=0$$ in $\R^N\setminus B(R)$. By local elliptic regularity
(\cite[Theorem 6.5.3]{krylov}), $u\in C^{2,\beta}_{loc}(\R^N)$ for every $\beta<1$. In
$\R^N\setminus B(R)$, $u$ satisfies
 $$\Delta u=\frac{\lambda u}{1+|x|^\alpha}\in L^p(\R^N).$$
By elliptic regularity, $u\in W^{2,p}(\R^N\setminus B(R))$. If
$p>\ds\frac{N}{2}$, we immediately deduce $u\in C_0(\R^N\setminus
B(R))$. Otherwise $u\in L^{p_1}(\R^N\setminus B(R))$ where
$\ds\frac{1}{p_1}=\ds\frac{1}{p}-\frac{2}{N}$ (with the usual
modification when $p=N/2)$. As before it follows $\Delta u\in
L^{p_1}(\R^N\setminus B(R))$ and $u\in W^{2,p_1}(\R^N\setminus
B(R))$. By iterating this procedure until $p_i>\frac{N}{2}$ we
deduce $u\in C_0(\R^N)$. Therefore $u$ attaints its minimum in a
point $x_0 \in \R^N$. The equality
$$\lambda u(x_0)=(1+|x_0|^\alpha)\Delta u(x_0)+f(x_0)$$
shows that $u(x_0) \ge 0$, since $\lambda >0$, hence $u \ge 0$. If
$\lambda =0 \in \rho(\hat{L}_p)$, then $\lambda \in \rho(\hat{L}_p)$
for small positive values of $\lambda$ and the thesis follows by
approximation. \qed

\begin{lem}  \label{confronto}
Suppose that $\lambda \in\rho(\hat{L}_p)$ for some $\lambda\geq 0$.
Then for every $0\leq f\in C_c(\R^N)$,
$$(\lambda-\hat{L}_p)^{-1}f\geq (\lambda-L_{min})^{-1}f.$$
\end{lem}
{\sc Proof.} Let $0\leq f\in C_c^\infty(\R^N)$  and set
$u=(\lambda-\hat{L}_p)^{-1}f$. Proposition \ref{positiv} and its
proof show that that $0\leq u\in D_{max}(L)$. Since
$(\lambda-L_{min})^{-1}f$ is the minimal solution, we immediately
have $u\geq (\lambda-L_{min})^{-1}f$.\qed

\begin{prop} \label{limitazione}
Let $N\geq 3$, $\alpha>2$, $p\leq \frac{N}{N-2}$. Then $\rho
(\hat{L}_p) \cap [0, \infty[=\emptyset$.
\end{prop}
{\sc Proof.} Let $\lambda>0$ and  $\chi_{B(0)}\leq f\leq\chi_{B(1)}$
be a smooth radial function. Denote by $u$ the (minimal) solution of
$\lambda u-Lu=f$ in $C_0(\R^N)$. Observe that from (\cite[Example
7.3]{met-wack}) it follows that the above equation has a unique
solution in $C_0(\R^N)$ (not in $C_b(\R^N)$). Hence, since the datum
$f$ is radial, the solution $u$ is radial too, and solves
$$\lambda u(\rho)-(1+\rho^\alpha)\left(u''(\rho)+\frac{N-1}{\rho}u'(\rho)\right)=f(\rho).$$
For $\rho\geq 1 $, $u$ solves the homogeneous equation
$$\lambda u(\rho)-(1+\rho^\alpha)\left(u''(\rho)+\frac{N-1}{\rho}u'(\rho)\right)=0.$$
Let us write $u$ as $u(\rho)=\eta(\rho)\rho^{2-N}$ for a suitable function $\eta$. Elementary computations show that $\eta$ satisfies
\begin{equation}   \label{ausiliaria}
\lambda \eta(\rho) -(1+\rho^\alpha)\left(\eta''(\rho)+\frac{3-N}{\rho}\eta'\right)=0
\end{equation}
for $\rho\geq 1$. First observe that, since $f\neq 0$ is
nonnegative, the strong maximum principle, see \cite[Theorem
3.5]{gil-tru}), implies that $u$ (and so $\eta$) is strictly
positive. We use Feller's theory to study the asymptotic behavior of
the solutions of the previous equation (see \cite[Section
VI.4.c]{engel-nagel}). We introduce the Wronskian
$$W(\rho)=\exp{\left\{-\int_1^\rho \frac{3-N}{s}ds\right\}}=\rho^{N-3}$$
and the functions
$$Q(\rho)=\frac{1}{(1+\rho^\alpha)W(\rho)}\int_1^\rho W(s)ds=\frac{1}{N-2}\frac{1}{(1+\rho^\alpha)\rho^{N-3}}(\rho^{N-2}-1)$$and
$$R(\rho)=W(\rho)\int_1^\rho\frac{1}{(1+s^\alpha)W(s)}ds=\rho^{N-3}\int_1^\rho\frac{1}{(1+s^\alpha)s^{N-3}}ds.$$
Since $\alpha>2$ by assumption, we have $Q\in L^1(1,+\infty)$ and
$R\notin L^1(1,\infty)$. This means that $\infty$ is an entrance
endpoint. In this case there exists a positive decreasing solution
$\eta_1$ of (\ref{ausiliaria}) satisfying
$\lim_{\rho\to\infty}\eta_1(\rho)=1$ and every solution of
(\ref{ausiliaria}) independent of $\eta_1$ is unbounded at
infinity. This shows that our solution $u$ grows at infinity at
least as $\rho^{2-N}$ and therefore it does not belong to
$L^p(\R^N)$. By Lemma \ref{confronto} we deduce that $\lambda \not
\in \rho(\hat{L}_p)$. \qed

When $N=1,2$ and $\alpha>2$, then (\ref{problem}) is never solvable
in $L^p$.

\begin{prop}
Let $N=1,2$, $\alpha>2$. Then $\rho (\hat{L}_p) \cap [0,
\infty[=\emptyset$.
\end{prop}
This follows from Proposition \ref{proprietaL} (ii), using Lemma
\ref{confronto}. \qed

\section{Solvability in $L^p$}
In this section we investigate the solvability of the equation
$\lambda u-Lu=f$ in $L^p$, for $\lambda \ge 0$. We start with
$\lambda=0$. Since the equation $-Lu=f$ is equivalent to $-\Delta
u(x)=f(x)/(1+|x|^\alpha)$, we can express $u$  and its gradient
through an integral operator involving the Newtonian potential. For
$f \in L^p$ we set
\begin{equation} \label{defu}
Tf(x)=u(x)=C_N\int_{\R^N} \frac{f(y)\,
dy}{(1+|y|^\alpha)|x-y|^{N-2}}
\end{equation}
and
\begin{equation} \label{defnablau}
Sf(x)=\nabla u(x)=C_N(N-2)\int_{\R^N} \frac{f(y) (y-x)\,
dy}{(1+|y|^\alpha)|x-y|^N}
\end{equation}
where $C_N=\left (N(2-N)\omega_N\right)^{-1}$ and $\omega_N$ is the
Lebesgue measure of the unit ball in $\R^N$.\\

We prove a preliminary result which will  be useful  to prove
estimates for the norm of the operator $T$ in $L^p$.

\begin{lem}  \label{sol-fond}
Let $2<\beta<N$. Then
$$\frac{1}{N(2-N)\omega_N}\int_{\R^N}\frac{dy}{|x-y|^{N-2}|y|^\beta}=\frac{1}{(2-\beta)(N-\beta)}|x|^{2-\beta}.$$
\end{lem}
{\sc Proof.} Set
$$u(x)=\frac{1}{N(2-N)\omega_N}\int_{\R^N}\frac{dy}{|x-y|^{N-2}|y|^\beta}.$$
By writing $x,\ y$  in spherical coordinates, $x=s\eta$, $y=r\omega$, with $\eta,\ \omega\in S_{N-1}$, $s,\ r\in [0, +\infty)$, the expression of $u$ becomes
\begin{eqnarray*}
 u(s\eta)&=&\frac{1}{N(2-N)\omega_N}\int_{S_{N-1}}d\omega\int_0^\infty\frac{r^{N-1} dr}{|s\eta-r\omega|^{N-2}|r|^\beta}\\&=&
\frac{1}{N(2-N)\omega_N}\int_{S_{N-1}}d\omega\int_0^\infty\frac{r^{N-1-\beta} dr}{s^{N-2}\left|\eta-\ds\frac{r}{s}\omega\right|^{N-2}}\\&=&
\frac{1}{N(2-N)\omega_N}s^{2-\beta}\int_{S_{N-1}}d\omega\int_0^\infty\frac{\xi^{N-1-\beta} d\xi}{\left|\eta-\xi\omega\right|^{N-2}}.
\end{eqnarray*}
By the rotational invariance of the integral,
\begin{eqnarray*}
 u(s\eta)&=&
\frac{1}{N(2-N)\omega_N}s^{2-\beta}\int_{S_{N-1}}d\omega\int_0^\infty\frac{\xi^{N-1-\beta} d\xi}{\left|e_1-\xi\omega\right|^{N-2}}\\&=&
\frac{1}{N(2-N)\omega_N}s^{2-\beta}\int_{\R^N}\frac{dy}{|e_1-y|^{N-2}|y|^\beta}.
\end{eqnarray*}
where $e_1$ is the unitary vector in the canonical basis of
$S_{N-1}$. Therefore $u(x)=C|x|^{2-\beta}$ with
$$C=\frac{1}{N(2-N)\omega_N}\int_{\R^N}\frac{dy}{|e_1-y|^{N-2}|y|^\beta}.$$
To compute  the constant $C$ we note that $u$ solves
$$\Delta u=\frac{1}{|x|^\beta}$$ or, in spherical coordinates,
$$
u''(\rho)+\frac{N-1}{\rho}u'(\rho)=\frac{1}{\rho^\beta}. $$

Inserting $u(\rho)=C\rho^{2-\beta}$. we get
$C=\ds\frac{1}{(2-\beta)(N-\beta)}$. \qed

In the following lemma we investigate the boundedness of the
operators $T,S$ in weighted $L^p$-spaces. Even though we need here
only the boundedness of $T$ in the unweighted $L^p$-space, we prove
the general result which will be of a central importance in the next
sections.

\begin{lem} \label{chius} Let $\alpha \ge 2$ and $N/(N-2) <p <\infty$.
For every $0 \le \beta, \gamma$ such that $\beta\leq \alpha-2$,
$\beta <\frac{N}{p'}-2$ and $\gamma\leq \alpha-1$, $\gamma <
\frac{N}{p'}-1$, there exists a positive constant $C$ such that for
any $f\in L^p$
\begin{eqnarray*}
& &\||\cdot|^\beta u\|_{L^p(\R^N)}\leq  C\|f\|_{L^p(\R^N)};\\
& &\||\cdot|^\gamma \nabla u\|_{L^p(\R^N)}\leq  C\|f\|_{L^p(\R^N)},
\end{eqnarray*}
where $u$ is defined in (\ref{defu}).
\end{lem}
{\sc Proof.}  Set $x=s\eta$, $y=\rho \omega$ with $s, \rho\in
[0,+\infty)$, $\eta \omega\in S_{N-1}$, then
\begin{eqnarray*}
u(s\eta)&=&\frac{1}{N(2-N)\omega_N}\int_{S_{N-1}}d\omega\int_0^\infty \frac{f(\rho\omega)\,
\rho^{N-1} \,
d\rho}{(1+\rho^\alpha)|s\eta-\rho\omega|^{N-2}}\\&=&\frac{1}{N(2-N)\omega_N}\int_{S_{N-1}}d\omega\int_0^\infty
\frac{s^2f(s\xi\omega)\, \xi^{N-1}\,
d\xi}{(1+(s\xi)^\alpha)|\eta-\xi\omega|^{N-2}}.
\end{eqnarray*}
We compute the $L^p$ norm of $|\cdot|^{\beta} u(\cdot)$. We start by
integrating with respect to $s$ the inequality above. We have, using
Minkowski inequality for integrals,
\begin{align*}
&\bigg(\int_0^\infty |u(s\eta)|^p\,  s^{\beta p+N-1}\,
ds\bigg)^\frac{1}{p}\leq\\ &\leq \frac{1}{N(N-2)\omega_N}\int_{S_{N-1}}d\omega\int_0^\infty
\frac{\xi^{N-1}\, d\xi}{|\eta-\xi\omega|^{N-2}}\left(\int_0^\infty
\frac{|f(s\xi\omega)|^p\,s^{N-1+2p+\beta
p}\,ds}{(1+s^\alpha\xi^\alpha)^p}\right)^{\frac{1}{p}}\\&=
\frac{1}{N(N-2)\omega_N}\int_{S_{N-1}}d\omega\int_0^\infty
\frac{\xi^{N-1}d\xi}{|\eta-\xi\omega|^{N-2}\xi^{\frac{N}{p}+\beta+2}}\left(\int_0^\infty
\frac{|f(v\omega)|^p}{(1+v^\alpha)^p}v^{N-1+2p+\beta p}\,
dv\right)^{\frac{1}{p}}.\\
\end{align*}
By recalling that  $\beta \le \alpha-2$ and since
$$
\frac{v^{2+\beta }}{1+v^\alpha} \le
\left(\frac{2+\beta}{\alpha-2+\beta}\right)^\frac{2+\beta}{\alpha}\frac{\alpha-2+\beta}{\alpha+2\beta},
$$
we obtain
\begin{align*}
\bigg(\int_0^\infty &|u(s\eta)|^p\,  s^{\beta p+N-1}\,
ds\bigg)^\frac{1}{p}\leq  \left(\frac{2+\beta}{\alpha-2+\beta}\right)^\frac{2+\beta}{\alpha}\frac{\alpha-2+\beta}{\alpha+2\beta}\frac{1}{N(N-2)\omega_N}\times\\ &\times\int_{S_{N-1}}d\omega\int_0^\infty
\frac{\xi^{N-1}d\xi}{|\eta-\xi\omega|^{N-2}\xi^{\frac{N}{p}+\beta+2}}\left(\int_0^\infty
|f(v\omega)|^p v^{N-1}\,
dv\right)^{\frac{1}{p}}.\\
\end{align*}
Let us observe that, by Lemma \ref{sol-fond} and the assumption
$\beta<\frac{N}{p'}-2$, we have
\begin{equation} \label{c}
\frac{1}{N(N-2)\omega_N}\int_{S_{N-1}}d\omega \int_0^\infty
\frac{\xi^{N-1}d\xi}{|\eta-\xi\omega|^{N-2}\xi^{\frac{N}{p}+\beta+2}}=\frac{p^2}{(N+\beta p)(Np-N-\beta p-2p)}.
\end{equation}
 By applying Jensen's inequality with
respect to probability measures
$$\ds\frac{\xi^{N-1}}{c|\eta-\xi\omega|^{N-2}\xi^{\frac{N}{p}+\beta+2}}\,d\xi
\, d\omega ,$$ where $c$ is the right-hand side in (\ref{c}), we
obtain
\begin{align*}
\int_0^\infty &|u(s\eta)|^p\,  s^{\beta p+N-1}\, ds  \leq
\left(\frac{2+\beta}{\alpha-2+\beta}\right)^\frac{(2+\beta)p}{\alpha}\left(\frac{\alpha-2+\beta}{\alpha+2\beta}\right)^p\frac{c^{p-1}}{N(N-2)\omega_N}\times\\&\times\int_{S_{N-1}}d\omega\int_0^\infty
\frac{\xi^{N-1}d\xi}{|\eta-\xi\omega|^{N-2}\xi^{\frac{N}{p}+\beta+2}}\int_0^\infty
|f(v\omega)|^p v^{N-1}\,
dv.\\
\end{align*}
By integrating with respect to $\eta$ on $S_{N-1}$, we obtain
\begin{align*}
\int_{\R^N} &|u(x)|^p\,  |x|^{\beta p}\, ds \leq \left(\frac{2+\beta}{\alpha-2+\beta}\right)^\frac{(2+\beta)p}{\alpha}\left(\frac{\alpha-2+\beta}{\alpha+2\beta}\right)^p\frac{c^{p-1}}{N(N-2)\omega_N}\times\\&\times
\int_{S_{N-1}}d\omega\int_{S_{N-1}}d\eta \int_0^\infty
\frac{\xi^{N-1}d\xi}{|\eta-\xi\omega|^{N-2}\xi^{\frac{N}{p}+\beta+2}}\int_0^\infty
|f(v\omega)|^p v^{N-1}\,
dv.\\
\end{align*}
A simple change of variables gives
\begin{align*}
\int_{S_{N-1}}d\eta \int_0^\infty
\frac{\xi^{N-1}d\xi}{|\eta-\xi\omega|^{N-2}\xi^{\frac{N}{p}+\beta+2}}&=\int_{S_{N-1}}d\eta
\int_0^\infty \frac{t^{N-1}dt}{|\eta
t-\omega|^{N-2}t^{\frac{N}{p'}-\beta}}\\&=\int_{\R^N}
\frac{dy}{|y-\omega|^{N-2}|y|^{\frac{N}{p'}-\beta}}.
\end{align*}
By applying Lemma \ref{sol-fond} again it follows that that
$$\int_{\R^N} |u(x)|^p\,  |x|^{\beta p}\,
dx\leq C^p \int_{\R^N} |f(x)|^p\, dx$$
with
\begin{equation}  \label{cost-prec}
C=\frac{p^2}{(N+\beta p)(Np-N-\beta p-2p)}\left(\frac{2+\beta}{\alpha-2+\beta}\right)^\frac{2+\beta}{\alpha}\frac{\alpha-2+\beta}{\alpha+2\beta}.
\end{equation}
The $L^p$ norm of $|\cdot|^{\gamma} \nabla u(\cdot)$ is estimated in
a similar way but we shall not be as precise as before concerning the
constants. By the representation formula,
$$|\nabla u(x)|\leq C\int_{\R^N} \frac{|f(y)|\, dy}{(1+|y|^\alpha)|x-y|^{N-1}}.$$
and hence
\begin{eqnarray*}
|\nabla u(s\eta)|&\leq &C\int_{S_{N-1}}d\omega\int_0^\infty
\frac{|f(\rho\omega)|\, \rho^{N-1} \,
d\rho}{(1+\rho^\alpha)|s\eta-\rho\omega|^{N-1}}\\&=
&C\int_{S_{N-1}}d\omega\int_0^\infty \frac{s|f(s\xi\omega)|\,
\xi^{N-1}\, d\xi}{(1+(s\xi)^\alpha)|\eta-\xi\omega|^{N-1}}.
\end{eqnarray*}
By Minkowski inequality and since $\gamma \le \alpha-1$,
\begin{align*}
\bigg(\int_0^\infty |\nabla u(s\eta)|^p\, & s^{\gamma p+N-1}\,
ds\bigg)^\frac{1}{p}\\ &\leq C\int_{S_{N-1}}d\omega\int_0^\infty
\frac{\xi^{N-1}\, d\xi}{|\eta-\xi\omega|^{N-1}}\left(\int_0^\infty
\frac{|f(s\xi\omega)|^p\,s^{N-1+p+\gamma
p}\,ds}{(1+s^\alpha\xi^\alpha)^p}\right)^{\frac{1}{p}}\\&=
C\int_{S_{N-1}}d\omega\int_0^\infty
\frac{\xi^{N-1}d\xi}{|\eta-\xi\omega|^{N-1}\xi^{\frac{N}{p}+\gamma+1}}\left(\int_0^\infty
\frac{|f(v\omega)|^p}{(1+v^\alpha)^p}v^{N-1+p+\gamma p}\,
dv\right)^{\frac{1}{p}}\\& \leq C\int_{S_{N-1}}d\omega\int_0^\infty
\frac{\xi^{N-1}d\xi}{|\eta-\xi\omega|^{N-1}\xi^{\frac{N}{p}+\gamma+1}}\left(\int_0^\infty
|f(v\omega)|^p v^{N-1}\,
dv\right)^{\frac{1}{p}}.\\
\end{align*}
As before, the assumption $\gamma<\frac{N}{p'}-1$ imples that the
integral
\begin{equation*}
\int_{S_{N-1}}d\omega \int_0^\infty
\frac{\xi^{N-1}d\xi}{|\eta-\xi\omega|^{N-1}\xi^{\frac{N}{p}+\gamma+1}}=\int_{\R^N}
\frac{dy}{|\eta-y|^{N-1}|y|^{\frac{N}{p}+\gamma+1}}
\end{equation*}
is finite and independent of $\eta \in S^{N-1}$ (by the rotational
invariance of the integrands). By applying Jensen's inequality we
obtain
\begin{align*}
\int_0^\infty |\nabla u(s\eta)|^p\, & s^{\gamma p+N-1}\, ds \\& \leq
C\int_{S_{N-1}}d\omega\int_0^\infty
\frac{\xi^{N-1}d\xi}{|\eta-\xi\omega|^{N-1}\xi^{\frac{N}{p}+\gamma+1}}\int_0^\infty
|f(v\omega)|^p v^{N-1}\,
dv.\\
\end{align*}
Integration with respect to $\eta$ on $S_{N-1}$ yields
\begin{align*}
\int_{\R^N} &|\nabla u(x)|^p\,  |x|^{\gamma p}\, ds\\ &\leq
C\int_{S_{N-1}}d\omega\int_{S_{N-1}}d\eta \int_0^\infty
\frac{\xi^{N-1}d\xi}{|\eta-\xi\omega|^{N-1}\xi^{\frac{N}{p}+\gamma+1}}\int_0^\infty
|f(v\omega)|^p v^{N-1}\,
dv.\\
\end{align*}
A simple change of variables gives
\begin{align*}
\int_{S_{N-1}}d\eta \int_0^\infty
\frac{\xi^{N-1}d\xi}{|\eta-\xi\omega|^{N-1}\xi^{\frac{N}{p}+\gamma+1}}&=\int_{S_{N-1}}d\eta
\int_0^\infty \frac{t^{N-1}dt}{|\eta
t-\omega|^{N-1}t^{N-\frac{N}{p}-\gamma}}\\&=\int_{\R^N}
\frac{dy}{|y-\omega|^{N-1}|y|^{N-\frac{N}{p}-\gamma}}.
\end{align*}
By the assumptions on $\gamma$, the last integral is convergent and
independent of $\omega$. It follows that
$$\int_{\R^N} |\nabla u(x)|^p\,  |x|^{\gamma p}\,
dx\leq C \int_{\R^N} |f(x)|^p\, dx$$ and the proof is complete. \qed

By the previous lemma,  the following estimate for the $L^p$-norm of
the operator $T$ immediately follows.

\begin{cor}  \label{stima-p}
Let $\alpha \ge 2$ and $N/(N-2) <p <\infty$. Then $$\|T\|_p\leq
{\left(\frac{2}{\alpha-2}\right)}^{\frac{2}{\alpha}}\frac{\alpha-2}{\alpha}\frac{p^2}{N(Np-N-2p)}.$$
\end{cor}
{\sc Proof.}
The estimate follows by setting  $\beta=0$ in (\ref{cost-prec}).\qed

\begin{os}\label{stimainfinito}{\rm The estimate with the  constant $C$ given by (\ref{cost-prec}) is stable as $p \to \infty$ only
if $\beta>0$. On the other hand, the operator $T$ is bounded also in
$L^\infty$ (and its norm will be computed later in Proposition
\ref{norm-inf}). It is possible to prove that the operator $T$ (with
$\beta=0$) is of weak-type $p-p$ with $p=N/(N-2)$ and then
interpolate between $N/(N-2)$ and $\infty$ to obtain stable
estimates for large $p$. The weak-type estimate is deduced as
follows. Write $Tf$ as the Riesz potential $I_2$ applied to the
function $f(x)/(1+|x|^\alpha)$ to get, using the classical estimate
of the Riesz potentials through the Hardy-Littlewood maximal
function $M$,
$$
|Tf(x)| \le C \left (M\left
(\frac{f(\cdot)}{1+|\cdot|^\alpha}\right)(x)\right
)^{1-2/N}\left\|\frac{f(\cdot)}{1+|\cdot|^\alpha}\right\|^{2/N}_1
$$
then Holder inequality to control the $L^1$-norms in terms of the
$L^{N/(N-2)}$-norm of $f$ and the weak $1-1$ estimate for $M$. Such
a proof works only for $\beta=0$ and gives constants depending on
those of the Marzinkiewicz interpolation theorem and of the
Hardy-Littlewood maximal function.}

\end{os}

We can now prove the invertibility of $L$ on $D_{p,max}(L)$, defined
in (\ref{massimale}).

\begin{prop} \label{invertibile}
Let $\alpha > 2$ and $N/(N-2) <p <\infty$. The operator $L$ is
closed and invertible on $D_{p, max}(L)$ and the inverse of $-L$ is
the operator $T$ defined in (\ref{defu}).
\end{prop}
{\sc Proof.} The closedness of $L$ on $D_{p, max}(L)$ follows from
local elliptic regularity. If $u \in D_{p, max}(L)$ satisfies
$Lu=0$, then $\Delta u=0$ and then $u=0$, since $u \in L^p$. This
shows the injectivity of $L$. Finally, let $f\in L^p$ and $f_n \in
C_c^\infty (\R^N)$ be such that $f_n \to f$ in $L^p$. Then $u_n=Tf_n
\to u=Tf$ in $L^p$, since $T$ is bounded (apply Lemma \ref{chius}
with $\beta=0$). By elementary potential theory
$$\Delta u_n(x)=\frac{f_n(x)}{1+|x|^\alpha}$$
hence $u_n \in D_{p, max}(L)$ and $Lu_n=f_n$. By the closedness of
$L$, $u \in D_{p, max}(L)$ and $Lu=f$.

\begin{teo} \label{risolvente}
Let $\alpha > 2$, $N/(N-2) <p <\infty$ and $\lambda \ge 0$. The
operator $\lambda -L$ is invertible on $D_{p, max}(L)$ and its
inverse is a positive operator. Moreover, if $f \in L^p \cap
C_0(\R^N)$, then $(\lambda-L)^{-1}f=(\lambda-L_{min})^{-1}f$.
\end{teo}
{\sc Proof.} Let $\rho$ be the resolvent set of $(L,D_{p, max}(L))$
and observe that the proposition above shows that $0 \in \rho$.
Lemma \ref{positiv} with $\hat{D}=D_{p,max}(L)$ shows that if $0\le
\lambda \in \rho$, than $(\lambda-L)^{-1} \ge 0$ and hence, by the
resolvent equation, $(\lambda-L)^{-1} \le (-L)^{-1}=T$ and therefore
\begin{equation} \label{stimapositiva}
\|(\lambda-L )^{-1}\| \le \|T\|
\end{equation}
where the norm above is the operator norm in $L^p$. Let
$E=[0,\infty[\cap \rho$. Then $E$ is non empty and open in
$[0,\infty[$, since $\rho$ is open, and closed since the operator
norm of $(\lambda-L)^{-1}$ is bounded in $E$. Then $E=[0,\infty[$.
To show the consistency of the resolvents we take $f \in C_c^\infty
(\R^N)$ and let $u=(\lambda-L)^{-1}f$. As in Lemma \ref{positiv} we
see that $u \in C_0(\R^N)\cap C^{2,\beta}_{loc}(\R^N)$ for any
$\beta<1$. If $supp\ f \subset \R^N\setminus B(R)$, then the
equation $Lu=\lambda u$ holds outside $B(R)$ and shows that $u$
belongs to $D_{max}(L)\cap C_0(\R^N)$, which is the domain of
$L_{min}$ in $C_0(\R^N)$, see Proposition \ref{proprietaL} (iii).
Therefore $(\lambda-L_{min})^{-1}f=u=(\lambda-L)^{-1}f$. By density,
this equality extends to all functions $f \in L^p\cap C_0(\R^N)$.
\qed

It is worth mentioning that the resolvents of $L$ in
$L^p$ and $L^q$ are consistent, provided that $p,q >N/(N-2)$. This
easily follows from above, together with a simple approximation
argument, since both resolvents are consistent with
$(\lambda-L_{min})^{-1}$. Observe also that estimate
(\ref{stimapositiva}) shows only that the resolvent is bounded on
$[0, \infty[$ and is not sufficient to apply the Hille-Yosida
theorem and prove results for parabolic problems. This will be done
in the next secion, under further restrictions on the admitted
values for $p$.

\section{Sectoriality in $L^p$}

We prove that, for $2<\alpha \le(N-2)(p-1)$ and $N/(N-2) <p
<\infty$, $(L, D_{p,max}(L))$ generates a strongly continuous
semigroup of positive contractions, analytic for $\alpha<(N-2)(p-1)$,
which coincides whith the minimal semigroup in $L^p \cap C_0(\R^N)$.

\begin{teo}
 Let $N\geq 3$, $p>N/(N-2)$, $2 < \alpha \leq (p-1)(N-2)$. Then $(L,D_{p,max}(L))$ generates a positive semigroup of contractions in $L^p$.
If $\alpha <(p-1)(N-2)$, the semigroup is also analytic.
\end{teo}
{\sc Proof.} Take $f\in L^p(\R^N)$, $\rho>0$, $\lambda \in \C$ and
consider the Dirichlet problem in $L^p(B(\rho))$
\begin{equation} \label{palla}
\left\{
\begin{array}{ll}
\lambda u-Lu=f & \textrm{in}\ B(\rho),\\
u=0   & \textrm{on}\ \partial B(\rho).
\end{array}\right.
\end{equation}
According to Theorem 9.15 in \cite{gil-tru}, for $\lambda>0$ there
exists a unique solution $u_\rho$ in $W^{2,p}(B(\rho))\cap
W_0^{1,p}(B(\rho))$. In order to show that the above problem is
solvable for complex values of $\lambda$ and to obtain estimates
independent of $\rho$, we show that $e^{\pm i\theta}L$ is
dissipative in $B(\rho)$ for $0\le  \theta  \le \theta_0$ and a
suitable $0<\theta_0 \le \pi/2$. Set
$u^\star=\ov{u}_\rho|u_\rho|^{p-2}$. Multiply  $Lu_\rho$  by
$u^\star$ and integrate over $B(\rho)$. The integration by parts is
straightforward when $p\geq 2$. For $1<p<2$, $|u_\rho|^{p-2}$
becomes singular near the zeros of $u_\rho$. It is possible to prove
the the integration by parts is allowed also in this case (see
\cite{met-spi}). Notice also that all boundary terms vanish since
$u_\rho=0$ at the boundary. So we get
\begin{align*}
\int_{B(\rho)}Lu_\rho\, &u^\star
dx=-\int_{B(\rho)}(1+|x|^\alpha)|u_\rho|^{p-4}|Re(\ov{u}_\rho\nabla
u_\rho)|^2dx-
\int_{B(\rho)}(1+|x|^\alpha)|u_\rho|^{p-4}|Im(\ov{u}_\rho\nabla
u_\rho)|^2dx\\&-\int_{B(\rho)}\ov{u}_\rho|u_\rho|^{p-2}\nabla(1+|x|^\alpha)\nabla
u\; dx- (p-2)\int_{B(\rho)}(1+|x|^\alpha)|u_\rho|^{p-4}
\ov{u}_\rho\nabla u_\rho Re(\ov{u}_\rho\nabla u_\rho)dx.
\end{align*}
By taking the real and imaginary part of the left and the right hand
side, we have
\begin{align*}
Re\bigg(\int_{B(\rho)}Lu_\rho\, &u^\star
dx\bigg)=-(p-1)\int_{B(\rho)}(1+|x|^\alpha)|u_\rho|^{p-4}|Re(\ov{u}_\rho\nabla
u_\rho)|^2dx\\&-
\int_{B(\rho)}(1+|x|^\alpha)|u_\rho|^{p-4}|Im(\ov{u}_\rho\nabla
u_\rho)|^2dx-\int_{B(\rho)}|u_\rho|^{p-2}\nabla(1+|x|^\alpha)Re(\ov{u}_\rho\nabla
u_\rho)\; dx;
\end{align*}
\begin{align*}
Im\bigg(\int_{B(\rho)}Lu_\rho\, u^\star
dx\bigg)=&-(p-2)\int_{B(\rho)}(1+|x|^\alpha)|u_\rho|^{p-4}Im(\ov{u}_\rho\nabla
u_\rho) Re(\ov{u}_\rho\nabla
u_\rho)dx\\&-\int_{B(\rho)}|u_\rho|^{p-2}\nabla(1+|x|^\alpha)Im(\ov{u}_\rho\nabla
u_\rho)\; dx.
\end{align*}
By Hardy's inequality as stated in Proposition \ref{hardy},
\begin{align*}
&\left|\int_{B(\rho)}|u_\rho|^{p-2}\nabla(1+|x|^\alpha)Re(\ov{u}_\rho\nabla
u_\rho)\; dx\right|\leq
\alpha\int_{B(\rho)}|u_\rho|^{p-2}|x|^{\alpha-1}|Re(\ov{u}_\rho\nabla
u_\rho)|dx\\\leq&
\alpha\left(\int_{B(\rho)}|u_\rho|^{p-4}|x|^{\alpha}|Re(\ov{u}_\rho\nabla
u_\rho)|^2
dx\right)^\frac{1}{2}\left(\int_{B(\rho)}|u_\rho|^p|x|^{\alpha-2}dx\right)^\frac{1}{2}\\\leq&
\frac{p\alpha}{\alpha-2+N}
\int_{B(\rho)}|u_\rho|^{p-4}(1+|x|^{\alpha})|Re(\ov{u}_\rho\nabla
u_\rho)|^2 dx.
\end{align*}
It follows that
\begin{align*}
-Re\bigg(\int_{B(\rho)}Lu_\rho\, &u^\star
dx\bigg)\geq\left(p-1-\frac{p\alpha}{\alpha-2+N}\right)\int_{B(\rho)}(1+|x|^\alpha)|u_\rho|^{p-4}|Re(\ov{u}_\rho\nabla
u_\rho)|^2dx\\&+
\int_{B(\rho)}(1+|x|^\alpha)|u_\rho|^{p-4}|Im(\ov{u}_\rho\nabla
u_\rho)|^2dx
\end{align*}
and
\begin{align*}
\bigg|Im &\bigg(\int_{B(\rho)}Lu_\rho\, u^\star
dx\bigg)\bigg|\\&\leq(p-2)\left(\int_{B(\rho)}|u_\rho|^{p-4}|x|^{\alpha}|Re(\ov{u}_\rho\nabla
u_\rho)|^2
dx\right)^\frac{1}{2}\left(\int_{B(\rho)}|u_\rho|^{p-4}|x|^{\alpha}|Im(\ov{u}_\rho\nabla
u_\rho)|^2
dx\right)^\frac{1}{2}\\&+\alpha\int_{B(\rho)}|u_\rho|^{p-2}|x|^{\alpha-1}|Im(\ov{u}_\rho\nabla
u_\rho)|\; dx\leq\\&
\leq\left(p-2+\frac{p\alpha}{\alpha-2+N}\right)\left(\int_{B(\rho)}|u_\rho|^{p-4}|x|^{\alpha}|Re(\ov{u}_\rho\nabla
u_\rho)|^2
dx\right)^\frac{1}{2}\\&\times\left(\int_{B(\rho)}|u_\rho|^{p-4}|x|^{\alpha}|Im(\ov{u}_\rho\nabla
u_\rho)|^2 dx\right)^\frac{1}{2}.
\end{align*}
Setting
\begin{align*}
B^2=\int_{\R^N}|u_\rho|^{p-4}|x|^{\alpha}|Re(\ov{u}_\rho\nabla u_\rho)|^2 dx,\\
C^2=\int_{\R^N}|u_\rho|^{p-4}|x|^{\alpha}|Im(\ov{u}_\rho\nabla
u_\rho)|^2 dx,
\end{align*}
we proved that
$$-Re\bigg(\int_{\R^N}Lu_\rho\, u^\star dx\bigg)\geq\left(p-1-\frac{p\alpha}{\alpha-2+N}\right)B^2+C^2$$
and
$$\left|Im\bigg(\int_{B(\rho)}Lu_\rho\, u^\star dx\bigg)\right|\leq \left(p-2+\frac{p\alpha}{\alpha-2+N}\right)BC.$$
Observe that $p-1-\ds\frac{p\alpha}{\alpha-2+N}$ is positive for
$\alpha<(N-2)(p-1)$. In this case it is possible to determine a
positive constant $l_\alpha$, independent of $\rho$, such that
$$\left(p-1-\frac{p\alpha}{\alpha-2+N}\right)B^2+C^2\geq l_\alpha\left(p-2+\frac{p\alpha}{\alpha-2+N}\right)BC$$
and, consequently,
$$\left|Im\bigg(\int_{B(\rho)}Lu_\rho\, u^\star dx\bigg)\right|\leq l_\alpha^{-1}\left\{-Re\bigg(\int_{B(\rho)}Lu_\rho\, u^\star dx\bigg)\right\}.$$
If $\tan\theta_\alpha=l_\alpha$, then $e^{\pm i\theta}L$ is
dissipative in $B(\rho)$ for $0 \le \theta \le \theta_\alpha$. The
previous computations give  also the dissipativity of $L$ if
$\alpha=(N-2)(p-1)$. Let us introduce the sector
$$
\Sigma_\theta=\{\lambda \in \C\setminus \{0\}: |Arg\,  \lambda| <
\pi/2+\theta\}.
$$
It follows from \cite[Theorem I.3.9]{pazy}, that problem
(\ref{palla}) has a unique solution for every $\lambda \in
\Sigma_\theta$ and $0 \le \theta <\theta_\alpha$ and that there
exists a constant $C_\theta$, independent of $\rho$, such that the
solution $u_\rho$ satisfies
\begin{equation} \label{stima}
\|u_\rho\|_{L^p (B(\rho))} \le
\frac{C_\theta}{|\lambda|}\|f\|_{L^p}.
\end{equation}
In the case $\alpha =(N-2)(p-1)$ the solutions $u_\rho$ exist for
$Re\, \lambda >0$ and satisfy the estimate
$$
\|u_\rho\|_{L^p(B(\rho))} \le \frac{1}{Re\, \lambda}\|f\|_{L^p}.
$$
Moreover, if $\lambda >0$ then $u_\rho\leq 0$ if $f\leq 0$ in
$B(\rho)$. In fact,  multiplying the equation
$$\lambda u_\rho-Lu_\rho=f$$
by $(u^+_\rho)^{p-1}$, integrating over $B(\rho)$ and proceeding as
before we obtain
$$\lambda\int_{B(\rho)}(u_\rho^+)^p\ dx\leq \int_{B(\rho)}f(u_\rho^+)^{p-1}\ dx\leq 0.$$ Therefore $u_\rho^+=0$ and $u_\rho\leq 0$.

Next we use weak compactness arguments to produce a function $u \in
D_{p,max}(L)$ satisfying $\lambda u-Lu=f$. For definiteness, we
consider the case $\alpha <(N-2)(p-1)$, the other one being simpler,
and fix $\lambda \in \Sigma_\theta$, with $0<\theta <\theta_\alpha$.

Let us fix a radius $r$ and apply the interior $L^p$ estimates
(\cite[Theorem 9.11]{gil-tru}) together with (\ref{stima}) to the
functions $u_\rho$ with $\rho >r+1$
$$\|u_\rho\|_{W^{2,p}(B(r))}\leq C_1[\|\lambda u_\rho-Lu_\rho\|_{L^p(B(r+1))}+\|u_\rho\|_{L^p(B(r+1))}] \le C_2\|f\|_{L^p}.$$
By weak compactness and a diagonal argument, we can find a sequence
$\rho_n \to \infty$ such that the functions $(u_{\rho_n})$ converge
weakly in $W^{2,p}_{loc}$ to a function $u$. Clearly $u$ satisfies
$\lambda u-Lu=f$ and, by (\ref{stima}),
\begin{equation}
\label{stima1}
 \|u\|_{L^p}\le
\frac{C_\theta}{|\lambda|}\|f\|_{L^p}.
\end{equation}
In particular $u \in D_{p,max}(L)$ and, moreover, $u$ is positive if
$\lambda, f \ge 0$. To complete the proof we need only to show that
$\lambda-L$ is injective on $D_{p,max}(L)$ for $\lambda \in
\Sigma_\theta$. Let
$$E=\{r>0: \Sigma_\theta \cap B(r) \subset\rho (L,D_{p,max}(L))\}$$
and $R=\sup E$. Since $0 \in E$, by Proposition \ref{invertibile},
$R$ is positive. On the other hand the norm of the resolvent exists
in $B(R) \cap \Sigma_\theta$ and is bounded by $C_\theta/|\lambda|$,
by (\ref{stima1}), hence cannot explode on the boundary of $B(R)$.
This proves that $R=\infty$ and concludes the proof. \qed

Finally, let us show that on $L^p \cap C_0(\R^N)$ the semigroup
coincide with $T_{min}$ of Section 2. In particular, the semigroups
are coherent in different $L^p$ spaces (when they are defined).

\begin{cor} \label{coerenza}
Let $(T_p(t))$ be the semigroup generated by $(L, D_{p,max}(L))$ in
$L^p$ and $(T_{min}(t)) $ be the minimal semigroup in $C_0(\R^N)$.
Then for every $f \in C_0(\R^N)\cap L^p$, $T_p(t)f=T_{min}(t)f$.
Moreover, if $p,q$ are allowed in the above theorem and $f \in
L^p\cap L^q$, then $T_p(t)f=T_q(t)f$.
\end{cor}
{\sc Proof.} Since  $(\lambda-L)^{-1}f=(\lambda-L_{min})^{-1}f$ for
$f \in L^p \cap C_0(\R^N)$, see Theorem \ref{risolvente}, the thesis
follows by representing the semigroups as the limit of iterates of
the corresponding resolvents. \qed

\section{Domain Characterization} The main result of this
section  consists in  showing that, for $N\geq 3$, $1<p<\infty$,
$2<\alpha<N/p'$,  the maximal domain $D_{p,max}(L)$ defined in
(\ref{massimale}) coincides with the weighted  Sobolev space $D_p$
defined by
$$
D_p=\{u\in W^{2,p}_{loc}(\R^N)\cap L^p(\R^N):\ (1+|x|^{\alpha-2})u,\
(1+|x|^{\alpha-1})\nabla u,\ (1+|x|^{\alpha}) D^2u\in L^p(\R^N)\}
$$
and endowed with its canonical norm.

\begin{os}
Observe that the assumption $2<\alpha<N/p'$ forces $p$ to be
strictly greater than $\frac{N}{N-2}$, according with Proposition
\ref{limitazione}.
\end{os}

\bigskip \noindent
The next lemma provides a core for $D_p$.
\begin{lem}  \label{densita}
The space $C_c^\infty(\R^N)$ is dense in $D_p$ with respect to
the norm
$$\|u\|_{D(A_p)}=\|u\|_{L^p(\R^N)}+\|(1+|x|^{\alpha-2})u\|_{L^p(\R^N)}+\|(1+|x|^{\alpha-1})\nabla u\|_{L^p(\R^N)}+\|(1+|x|^{\alpha})D^2 u\|_{L^p(\R^N)}.$$
\end{lem}
{\sc Proof.} Let us first observe that a function $u\in
W^{2,p}(\R^N)$ with compact support can be approximated by a
sequence of $C^\infty$ functions with compact support, in the
$D(A_p)$ norm. Indeed, if $\rho_n$ are standard mollifiers,
$u_n=\rho_n\ast u\in C_c^\infty(\R^N)$, $\supp u_n\subset \supp
u+B(1)$ for any $n\in N$ and  $u_n\to u$ in $D_p$ since
$(1+|x|^{\alpha-2}),\ (1+|x|^{\alpha-1}),\ (1+|x|^\alpha)$ are
bounded (uniformly with respect to $n$) on $\supp u+B(1)$. Next we
show that any function $u$ in $D_p$ can be approximated, with
respect to the norm of $D(A_p)$, by a sequence of functions in
$W^{2,p}(\R^N)$ each having a compact support. Let $\eta$ be a
smooth function such that $\eta=1$ in $B(1)$, $\eta=0$ in
$\R^N\setminus B(2)$, $0\leq \eta\leq 1$ and set
$\eta_n(x)=\eta\left(\frac{x}{n}\right)$. If $u\in D_p$, then
$u_n=\eta_n u$ are compactly supported functions in $W^{2,p}(\R^N)$,
$u_n\to u$ in $L^p(\R^N)$, $(1+|x|^{\alpha-2})u_n\to
(1+|x|^{\alpha-2})u$ in $L^p(\R^N)$ by dominated convergence.
Concerning the convergence of the derivatives we have
$$(1+|x|^{\alpha-1})\nabla u_n=\frac{1}{n}(1+|x|^{\alpha-1})\nabla\eta\left(\frac{x}{n}\right)u+(1+|x|^{\alpha-1})\eta\left(\frac{x}{n}\right)\nabla u.$$
As before, $$(1+|x|^{\alpha-1})\eta\left(\frac{x}{n}\right)\nabla u
\to (1+|x|^{\alpha-1})\nabla u$$ in $L^p(\R^N)$. For the left term,
since $\nabla \eta (x/n)$ can be different from zero only for $n \le
|x| \le 2n$ we have
$$\frac{1}{n}(1+|x|^{\alpha-1})\left|\nabla\eta\left(\frac{x}{n}\right)\right||u|\leq C(1+|x|^{\alpha-2})|u|\chi_{\{n\leq |x|\leq 2n\}},$$
and the right hand side tends to $0$ as $n\to\infty$. A similar
argument shows the convergence of the second order derivatives in
the weighted  $L^p$ norm.\qed

We can  prove that $L$ is closed on $D_p$.

\begin{prop} \label{chiusura}
assume that $2<\alpha <N/p'$. Then there exists a positive constant
$C$ such that for any $u\in D_p$
\begin{eqnarray*}
& &\|u\|_{L^p(\R^N)}+\|(1+|x|^{\alpha-2})u\|_{L^p(\R^N)}+\|(1+|x|^{\alpha-1})\nabla u\|_{L^p(\R^N)}+\|(1+|x|^{\alpha})D^2 u\|_{L^p(\R^N)} \\&\leq & C\|Lu\|_{L^p(\R^N)}.
\end{eqnarray*}

\end{prop}
{\sc Proof.} Let $u \in C_c^\infty (\R^N)$  and set $f=-Lu$. Then $f
\in C^2_c(\R^N)$, $(1+|x|^\alpha)\Delta u(x)=-f(x)$ or equivalently
$$
-\Delta u (x) =\frac{f(x)}{1+|x|^\alpha}.
$$
By elementary potential theory, $u$ is given by (\ref{defu}). By
setting $\beta=\alpha-2$ and $\gamma=\alpha-1$ in Lemma \ref{chius}
and since $\alpha<\frac{N}{p'}$, we deduce that
$$\|u\|_{L^p(\R^N)}\leq C\|Lu\|_{L^p(\R^N)},$$

$$\|(1+|x|^{\alpha-2})u\|_{L^p(\R^N)}\leq C\|Lu\|_{L^p(\R^N)},$$
and
$$\|(1+|x|^{\alpha-1})\nabla u\|_{L^p(\R^N)}\leq C\|Lu\|_{L^p(\R^N)}$$
for every $u \in C_c^\infty (\R^N)$. In order to prove the estimates
of the second order derivatives, we  apply the classical
Calder\'{o}n- Zygmund inequality to $(1+|x|^\alpha)u$ and the
estimates of the lower order derivates obtained above. We deduce
\begin{align*}
\|(1+|x|^\alpha)&D^2 u\|_{L^p(\R^N)}\\&\leq
C(\alpha)[\|D^2((1+|x|^\alpha)u)\|_{L^p(\R^N)})+\|(1+|x|^{\alpha-1})\nabla
u\|_{L^p(\R^N)} +\|(1+|x|^{\alpha-2})u\|_{L^p(\R^N)}]\\&\leq
C(N,p,\alpha)[\|\Delta((1+|x|^\alpha)
u)\|_{L^p(\R^N)}+\|Lu\|_{L^p(\R^N)}]\\&\leq
C(N,p,\alpha)[\|(1+|x|^\alpha) \Delta
u\|_{L^p(\R^N)}+\|Lu\|_{L^p(\R^N)}]=C(N,p,\alpha)\|Lu\|_{L^p(\R^N)}.
\end{align*}
By Lemma \ref{densita}, these estimates extend from
$C_c^\infty(\R^N)$ to $D_p$. \qed

The following lemma is a tool to prove the equality $D_p=D_{p,max}(L)$.
Once the latter equality has been proved, it is an obvious
consequence of Proposition \ref{invertibile}.
\begin{lem}  \label{invertib}
If $2<\alpha <N/p'$, the operator $-(L,D_p)$ is invertible and its
inverse is the operator $T$ defined in (\ref{defu}).
\end{lem}
{\sc Proof.} In fact, $T$ is bounded in $L^p$, by Lemma \ref{chius}
with $\beta=0$, and the equality $u=-TLu$ holds for every $u \in
C_c^\infty (\R^N)$. Since $C_c^\infty (\R^N)$ is a core for
$(L,D_p)$, see Lemma \ref{densita}, then $u=-TLu$ for every $u \in
D_p$. Since $(L,D_p)$ is injective, the proof is complete. \qed

\begin{teo} \label{dominio}
If $2<\alpha <N/p'$, then $D_p$ coincides with the maximal domain in
$L^p$, that is
$$D_p=\{u\in L^p\cap W^{2,p}_{loc}:\ Lu\in L^p\}.$$
\end{teo}
{\sc Proof.} The inclusion $D_p\subset D_{p,max}(L)$ is obvious. Let
now $u\in D_{p,max}(L)$. By Corollary \ref{invertib}, there exists
$v\in D_p$ such that $Lv=Lu$. Therefore $u-v$ belongs to the maximal
domain of $L$ and $L(u-v)=0$, that is $\Delta (u-v)=0$ . Since $u,v
\in L^p$, then $u=v$ and $u$ belongs to $D_p$.\qed

Next we show that if $\alpha \ge N/p'$ then $D_p$ is properly
contained in $D_{p,max}(L)$

\begin{prop} Let $N\geq 3$, $p>N/(N-2)$, $\alpha\geq N/p'$, $\alpha >2$. Then
$D_p$ is a proper subset of $D_{p,max}(L)$.
\end{prop}
{\sc Proof.} Let $\chi_{B(1)}\leq f\leq\chi_{B(2)}$ be a smooth
radial function. Denote by $u$ the unique solution in $D_{p,max}(L)$
of $Lu(\rho)=(1+\rho^\alpha)f(\rho)$, see Proposition
\ref{invertibile}. Since the datum $f$ is radial, by uniqueness, the
solution $u$ is radial too, hence it solves
$$u''(\rho)+\frac{N-1}{\rho}u'(\rho)=f(\rho).$$
For $\rho\geq 2 $, $u$ solves the homogeneous equation
$$u''(\rho)+\frac{N-1}{\rho}u'(\rho)=0,$$ hence it is given by $u=c\rho^{2-N}$ for some positive $c$. Then
$$\int_{\R^N}|(1+|x|^{\alpha-2})u|^p dx\geq c_1\int_2^\infty (1+\rho^{\alpha-2})^p\rho^{p(2-N)}\rho^{N-1}
\ d\rho\geq C\int_2^\infty \rho^{p\alpha-Np+N-1}d\rho.$$ The last
integral converges if and only if $\ds\alpha<\frac{N}{p'}$. In a
similar way one can show that  $(1+|x|^{\alpha-1})\nabla u$ and
$(1+|x|^{\alpha}) D^2u$ are not in $L^p(\R^N)$. \qed

A partial characterization of $D_{p,max}(L)$ can be obtained from
Lemma \ref{chius}.
\begin{prop} \label{partial}
Let $N\geq 3$, $p>N/(N-2)$, $\alpha\geq N/p'$. If $0 \le \beta
<N/p'-2$ and $0 \le \gamma <N/p'-1$, then $|x|^\beta u$ and
$|x|^\gamma \nabla u$ belong to $L^p$, for every $u \in
D_{p,max}(L)$.
\end{prop}
{\sc Proof.} This follows immediately from Lemma \ref{chius}, since
the operator $-T$ defined in (\ref{defu}) is the inverse of
$(L,D_{p,max}(L))$. \qed

Observe that for $\beta=\gamma=0$ the above result has been already
proved in Proposition \ref{w2p}.

If $\alpha >N/p'$, then $C_c^\infty (\R^N)$ is not a core for $L$.
This fact also gives $D_p \neq D_{p,max}(L)$ in this case.

\begin{prop}
Let $\alpha>N/p'$. Then $L(C_c^\infty(\R^N))$ is not dense in
$L^p(\R^N)$.
\end{prop}
{\sc Proof.} It is sufficient to observe that $0\neq \ds\frac{1}{1+|x|^\alpha}\in L^{p'}(\R^N)$ and
$$\int_{\R^N} (1+|x|^\alpha)\Delta u \frac{1}{1+|x|^\alpha}\ dx=\int_{\R^N}\Delta u\ dx=0$$
for every $u\in C_c^\infty(\R^N)$.\qed

\begin{prop}   \label{dens}
Let $\alpha=N/p'$. Then $L(C_c^\infty(\R^N))$ is dense in $L^p(\R^N)$.
\end{prop}
{\sc Proof.} Let $g\in L^{p'}(\R^N)$ such that
$$\int_{\R^N}(1+|x|^\alpha)\Delta u\cdot g\, dx=0$$
for every $u\in C_c^\infty(\R^N)$. It follows that
$$\Delta (g(1+|x|^\alpha))=0$$
in the distributional sense, hence in a classical sense. Set
$h=g(1+|x|^\alpha)$. Since $h$ is an harmonic function,   it
satisfies
$$|\nabla h(0)|\leq \frac{C}{R^{N+1}}\int_{B(0, R)}|h|dx $$
for every  $R>0$. By assumption $g=\ds\frac{h}{(1+|x|^\alpha)}\in
L^{p'}(\R^N)$, therefore H\"{o}lder's inequality yields
 $$|\nabla h(0)|\leq \frac{C}{R^{N+1}}\int_{B(0, R)}\frac{|h|}{1+|x|^\alpha}(1+|x|^\alpha)dx
  \leq CR^{-N-1+\frac{N}{p}+\alpha}=CR^{\alpha-\frac{N}{p'}-1}=CR^{-1}.$$
Letting $R$ to infinity, we deduce that $|\nabla h(0)|=0$. In a
similar way one proves that $|\nabla h(x_0)|=0$ for every
$x_0\in\R^N$. It means that $h=C$ for some constant $C$ and
$g=\ds\frac{C}{1+|x|^\alpha}\in L^{p'}(\R^N)$. But
$\ds\frac{1}{1+|x|^\alpha}\in L^{p'}(\R^N)$ if and only if
$\alpha>\frac{N}{p'}$, therefore $C=0$ and, consequently, $g=0$.
This proves the density of  $L(C_c^\infty(\R^N))$ in $L^p(\R^N)$.
\qed

It can be proved that the a-priori estimates of Proposition
\ref{chiusura} for $p=2$ still hold in $C_c^\infty (\R^N)$ if
$\alpha \neq N/2$. However, since $C_c^\infty (\R^N)$ is not a core
for $(L,D_{2,max}(L)$, for $\alpha >N/2$ they do not extend to the
domain of $L$. Next we show that the a-priori estimates fail even in
$C_c^\infty (\R^N)$ if $\alpha =N/p'$, which is a core by the
Proposition above.

\begin{prop}
Let $N\geq 3,\ \alpha=N/p'$. Then the estimates in
Proposition \ref{chiusura} do not hold on $C_c^\infty (\R^N)$.
\end{prop}
{\sc Proof.} Let, $R\geq 2$, $\phi\in C_c^\infty(\R^N)$ be a radial
function such that $\phi=1$ in $B(R)\setminus B(2)$, $\phi=0$ in
$B(1) \cup \left (\R^N\setminus B(2R)\right )$,
$\|\phi'_R\|_\infty\leq\frac{C}{R}$,
$\|\phi''_R\|_\infty\leq\frac{C}{R^2}$ and set
$u(\rho)=\phi_R\rho^{2-N}$, $N>2$. Then $u \in C_c^\infty(\R^N)$ (we
omit to indicate the dependence of $u$ on $R$) and
$$\Delta
u=u''+\frac{N-1}{\rho}u'=\phi''_R(\rho)\rho^{2-N}+(3-N)\phi'_R(\rho)\rho^{1-N}.$$
A straightforward computation shows that, for $\alpha=\frac{N}{p'}$,
$$ \int_{\R^N}|(1+|x|^\alpha)\Delta u|^p\ dx \leq C
$$
with $C$ independent of $R$.  On the other hand
$u'(\rho)=\phi_R'(\rho)\rho^{2-N}+\phi_R(2-N)\rho^{1-N}$ and
$$\int_0^\infty (1+\rho^{\alpha-1})^p\rho^{N-1}|u'(\rho)|^p\ d\rho
=\int_1^{2R}
(1+\rho^{\alpha-1})^p\rho^{N-1}\left|\phi_R'(\rho)\rho^{2-N}+\phi_R(2-N)\rho^{1-N}\right|^p
d\rho.$$ The last integral tends to $\infty$ as $R \to \infty$ since
$$\int_1^{2R}
(1+\rho^{\alpha-1})^p\rho^{N-1}\left|\phi_R(2-N)\rho^{1-N}\right|^p
\, d\rho \ge C \log R.
$$ Therefore the $L^p$-norm of $(1+|x|^{\alpha-1})\nabla u$ cannot
be controlled by the $L^p$-norm of $Lu$ on $C_c^\infty(\R^N)$
Similarly one shows that the $L^p$-norm of $(1+|x|^{\alpha})D^2 u$
cannot be controlled by the $L^p$-norm of $Lu$.

\qed

\section{The operator in $C_0(\R^N)$}
Let
$$
D(L)= D_{max}(L)\cap C_0(\R^N)
$$
be the generator of $(T_{min}(t))_{t \ge 0}$ in $C_0$, see
Proposition \ref{proprietaL} (iii) and note that we need the only
restriction $N, \alpha >2$.\\

As in the $L^p$-case we give a description of the domain when
$\alpha <N$ and a partial description when $\alpha \ge N$.

We need the analogous of Lemma \ref{chius} for $p=\infty$

\begin{lem} \label{comportamentoJ}
Let $\gamma,\ \beta>0$ such that $\gamma<N$ and $\gamma+\beta>N$.
Set
$$J(x)=\int_{\R^N}\frac{dy}{|x-y|^\gamma(1+|y|^\beta)}.$$
 Then $J$ is bounded in $\R^N$ and has the following  behaviour as $|x|$ goes to infinity
\begin{equation*}
J(x)\simeq\left\{
\begin{array}{ll}
c_1 |x|^{N-(\gamma+\beta)}& \textrm{if}\quad\beta<N \\
c_2|x|^{-\gamma}\log|x|   &\textrm{if}\quad\beta=N\\
c_3|x|^{-\gamma}  &\textrm{if}\quad\beta>N \end{array}\right.
\end{equation*}
for suitable positive constants $c_1, c_2, c_3$.
 \end{lem}
{\sc Proof.} Since $\ds\frac{1}{|y|^{\gamma}}$ and
$\ds\frac{1}{1+|y|^\beta}$ are radial decreasing $J(x) \le J(0)
<\infty$. In order to prove the asymptotic behaviour, we write $J$
in spherical coordinates. Set $x=s\eta$, $y=\rho \omega$ with $s,
\rho\in [0,+\infty)$, $\eta \omega\in S_{N-1}$, then
\begin{eqnarray*}
J(s\eta)&=&\int_{S_{N-1}}d\omega\int_0^\infty \frac{\rho^{N-1} \,
d\rho}{(1+\rho^\alpha)|s\eta-\rho\omega|^{N-2}}\\&=&\int_{S_{N-1}}d\omega\int_0^\infty
\frac{s^N \xi^{N-1}\,
d\xi}{s^\gamma(1+(s\xi)^\beta)|\eta-\xi\omega|^{\gamma}}=\int_{S_{N-1}}d\omega\int_0^\infty
\frac{s^{N-\gamma} \xi^{N-1}\,
d\xi}{(1+s^\beta\xi^\beta)|e_1-\xi\omega|^{\gamma}}\\&=&\int_{S_{N-1}}d\omega\int_0^\frac{1}{2}
\frac{s^{N-\gamma} \xi^{N-1}\,
d\xi}{(1+s^\beta\xi^\beta)|e_1-\xi\omega|^{\gamma}}+\int_{S_{N-1}}d\omega\int_{\frac{1}{2}}^\infty
\frac{s^{N-\gamma} \xi^{N-1}\,
d\xi}{(1+s^\beta\xi^\beta)|e_1-\xi\omega|^{\gamma}}.
\end{eqnarray*}
Set $$J_1(s\eta)=\int_{S_{N-1}}d\omega\int_0^\frac{1}{2}
\frac{s^{N-\gamma} \xi^{N-1}\,
d\xi}{(1+s^\beta\xi^\beta)|e_1-\xi\omega|^{\gamma}}$$ and $$J_2(s\eta)=\int_{S_{N-1}}d\omega\int_{\frac{1}{2}}^\infty
\frac{s^{N-\gamma} \xi^{N-1}\,
d\xi}{(1+s^\beta\xi^\beta)|e_1-\xi\omega|^{\gamma}}.$$
Concerning $J_2$, we have
$$\lim_{s\to +\infty}s^{\gamma+\beta-N}J_2=\int_{S_{N-1}}d\omega\int_{\frac{1}{2}}^\infty
\frac{\xi^{N-1}\,
d\xi}{\xi^\beta|e_1-\xi\omega|^{\gamma}}=\int_{\R^N\setminus B(0,\frac{1}{2})}\frac{dy}{|y|^\beta|e_1-y|^\gamma}=C>0$$
for some positive contant $C$. Therefore
\begin{equation}  \label{asint1}
J_2(x)\simeq C|x|^{N-(\gamma+\beta)}
 \end{equation}
as $|x|\to\infty$. Let us estimate the remaining term. We have
\begin{eqnarray*}J_1(s\eta)=\int_{S_{N-1}}d\omega\int_0^\frac{1}{2}
\frac{s^{N-\gamma} \xi^{N-1}\,
d\xi}{(1+s^\beta\xi^\beta)|e_1-\xi\omega|^{\gamma}}&=&s^{-\gamma} \int_{S_{N-1}}d\omega\int_0^\frac{s}{2}
\frac{t^{N-1}\,
dt}{(1+t^\beta)\left|e_1-\frac{t}{s}\omega\right|^{\gamma}}.
\end{eqnarray*}
Since $\frac{1}{2}\leq\left|e_1-\frac{t}{s}\omega\right|\leq \frac{3}{2}$,
$$c_1 J_1\leq s^{-\gamma}\int_0^\frac{s}{2}
\frac{t^{N-1}\,
dt}{(1+t^\beta)}\leq c_2 J_1$$
for some positive $c_1,\ c_2$.
Evidently
\begin{equation*}
s^{-\gamma}\int_0^\frac{s}{2}
\frac{t^{N-1}\,
dt}{(1+t^\beta)}\simeq\left\{
\begin{array}{ll}
 |s|^{N-(\gamma+\beta)}& \textrm{if}\quad\beta<N, \\
|s|^{-\gamma}\log|x|   &\textrm{if}\quad\beta=N,\\
|s|^{-\gamma}  &\textrm{if}\quad\beta>N \end{array}\right.
\end{equation*}
as $s$ goes to infinity. From (\ref{asint1}) and the last estimate
the aymptotic behaviour of $J$ follows.
\qed

The following two results are deduced from the lemma above as
Theorem \ref{dominio} and Proposition \ref{partial} are deduced from
Lemma \ref{chius}.

\begin{teo} \label{dominio-infinito}
Let $2<\alpha<N$. Then
$$
D(L)=\{u\in C_0:\ (1+|x|^{\alpha-2})u,\ (1+|x|^{\alpha-1})\nabla u,\
(1+|x|^{\alpha}) \Delta u\in C_0\}
$$
\end{teo}

\begin{prop} \label{partial1}
Let $N\geq 3$,  $\alpha\geq N$. If $0 \le \beta <N-2$ and $0 \le
\gamma <N-1$, then for every $u \in D(L)$, $|x|^\beta u$ and
$|x|^\gamma \nabla u$ belong to $C_0$.
\end{prop}

Finally, we compute the operator norm in $C_0$ of the operator
$T=(-L)^{-1}$ defined in (\ref{defu})
\begin{prop} \label{norm-inf}
If $N \ge 3$ and $\alpha >2$ then
\begin{equation*}
\|T\|_\infty= \frac{\pi}{(N-2)\alpha\sin\left(\frac{2}{\alpha}\pi\right)}.
\end{equation*}
\end{prop}
{\sc Proof.} We have
$$
\|T\|=\frac{1}{N(N-2)\omega_N}\sup_{x \in
\R^N}J(x)=\frac{1}{N(N-2)\omega_N}J(0)=
\frac{1}{N(N-2)\omega_N}\int_{\R^N}\frac{dy}{|y|^{N-2}(1+|y|^\alpha)}.$$
Since
$$\frac{1}{N\omega_N}\int_{\R^N}\frac{dy}{|y|^{N-2}(1+|y|^\alpha)}=\int_0^\infty\frac{s^{N-1}}{s^{N-2}(1+s^\alpha)}ds=
\frac{1}{\alpha}\int_0^\infty\frac{t^{\frac{2}{\alpha}-1}}{1+t}dt=\frac{\pi}{\alpha
\sin\left(\frac{2}{\alpha}\pi\right)}$$ the proof is complete. \qed

\section{Discreteness and location of the spectrum}
Throughout this section, to unify the notation, when $p=\infty$,
$L^p$ stands for $C_0$ and $D_{p,max}(L)$ for $D(L)$.
\begin{prop} \label{comp}
If $N/(N-2)<p<\infty$, $2<\alpha \le \infty$, then the resolvent of
$(L,D_{p,max}(L))$ is compact in $L^p$.
\end{prop}
{\sc Proof.} Let us prove that $D_{p,max(L)}$ is compactly embedded
into $L^p$ for $p<\infty$. Let $\mathcal{U}$ be a bounded subset of
$D_{p,max}(L)$. Fixing $0<\beta <\alpha-2,\ N/p'-2$ in Lemma
\ref{chius} we obtain $\int_{\R^N}|(1+|x|^{\beta})u|^p\leq M$ for
some positive $M$ and for every $u\in\mathcal{U}$. Then, given
$\eps>0$, there exists $R>0$ such that
$$\int_{|x|>R}|u|^p<\eps^p$$
for every $u\in \mathcal{U}$. Let $\mathcal{U}'$ be the set of the restrictions of the functions in $\mathcal{U}$ to $B(R)$.
Since the embedding of $W^{2,p}(B(R))$ into $L^p(B(R))$ is compact, the set $\mathcal{U}'$ which is bounded in $W^{2,p}(B(R))$ is totally bounded in $L^p(B(R))$. Therefore there exist $n\in\N$, $f_1,\ldots,f_n\in L^p(B(R))$ such that
$$\mathcal{U}'\subseteq \bigcup_{i=1}^n\{f\in L^p(B(R)):\ \|f-f_i\|_{L^p(B(R))}<\eps\}.$$
Set $\ov{f}_i=f_i$ in $B(R)$ and $\ov{f}_i=0$ in $\R^N\setminus B(R)$. Then $\ov{f}_i\in L^p(\R^N)$ and
$$\mathcal{U}\subseteq \bigcup_{i=1}^n\{f\in L^p(\R^N):\ \|f-\ov{f}_i\|_{L^p(\R^)}<2\eps\}.$$
It follows that $\mathcal{U}$ is relatively compact in $L^p(\R^N)$.
The compactness of the resolvent of $(L,D(L))$ in $C_0$ follows
similarly from the results of the previous section or from
(\cite[Example 7.3]{met-wack}). \qed

Clearly, the spectrum of $L$ consists of eigenvalues. Let us show
that it is independent of $p$.

\begin{cor} \label{indipendenza}
If $N/(N-2)<p \le \infty$, $2<\alpha<\infty$, then the spectrum of
$(L,D_{p,max}(L))$ is independent of $p$.
\end{cor}
{\sc Proof.}  Let $\rho_p, \rho_q$ be the resolvent sets in $L^p, L^q$, respectively.
Then $0 \in \rho_p \cap \rho_q$ and the inverse of $L$ in $L^p$ and in $L^q$
 is given by the operator $-T$ defined in (\ref{defu}), see Proposition \ref{invertibile}.
 This shows the consistency of the resolvents at $0$ and, since $\rho_p\cap \rho_q$ is connected,
 the consistency of the resolvents at any point of $\rho_p\cap \rho_q$, see \cite[Proposition 2.2]{arendt}.
 An application of \cite[Proposition 2.6]{arendt}  concludes the
 proof.
\qed
 In order to have more information on the spectrum of $L$, we
 introduce the Hilbert space $L^2_\mu$, where
 $d\mu(x)=(1+|x|^\alpha)^{-1}dx$, endowed with its canonical inner
 product. Note that the measure $\mu$ is finite if and only if
 $\alpha >N$. We consider also the Sobolev space
 $$
 H=\{u \in L^2_\mu : \nabla u \in L^2\}
 $$
 endowed with the inner product
 $$
 (u,v)_H=\int_{\R^N}\left (u\bar{v}\, d\mu +\nabla u \cdot \nabla
 \bar{v}\, dx \right)
 $$
and let $\cal V$ be the closure of $C_c^\infty$ in $H$, with respect
to the norm of $H$. Observe that Sobolev inequality
\begin{equation} \label{sobolev}
\|u\|_{2^*}^2 \le C_2^2 \|\nabla u\|_2^2
\end{equation}
holds in $\cal V$ but not in $H$ (consider for example the case
where $\alpha >N$ and $u=1$). Here $2^*=2N/(N-2)$ and $C_2$ is the
best constant for which the equality above holds.

\begin{lem} \label{compatto}
If $\alpha >N$, the embedding of $\cal V$ in $L^2_\mu$ is compact.
\end{lem}
{\sc Proof.} The proof is very similar to that of Proposition
\ref{comp} once one notes that on any ball $B(R)$ the measure $\mu$
is bounded above and below from zero. Therefore, it suffices to show
that given $\mathcal{U}$ a bounded subset of $\cal V$ and $\eps >0$,
there exists $R>0$ such that
$$\int_{|x|>R}|u|^2\, d\mu<\eps^2$$
for every $u\in \mathcal{U}$. This easily follows from
(\ref{sobolev}) since
$$
\int_{|x|>R}|u|^2\, d\mu \le \left
(\int_{|x|>R}|u|^{\frac{2N}{N-2}}\,dx \right )^{1-\frac{2}{N}}\left
(\int_{|x|>R}\frac{1}{(1+|x|^\alpha)^{\frac{N}{2}}}\,dx \right
)^{\frac{2}{N}}.
$$
\qed

Next we introduce the continuous and weakly coercive symmetric form
\begin{equation} \label{forma}
a(u,v)= \int_{\R^N}\nabla u \cdot \nabla
 \bar{v}\, dx
\end{equation}
for $u,v \in \cal V$ and the self-adjoint operator $\cal L$
 defined by
 $$
 D({\cal L})=\{u\in L^2_\mu : {\rm there\ exists\ } f \in L^2_\mu:
 a(u,v)=-\int_{\R^N}f\bar{v}\, d\mu {\rm \ for\ every\  } v \in {\cal V} \} \qquad {\cal
 L}u=f.
$$
Since $a(u,u) \ge 0$, the operator ${\cal L}$ generates an analytic
semigroup of contractions in $e^{t{\cal L}} $ in $L^2_\mu$. An
application of the Beurling-Deny criteria shows that the generated
semigroup is positive and $L^\infty$-contractive. For our purposes
we need to show that the resolvent of ${\cal L}$ and of
$(L,D_{p,max}(L))$  are coherent. This is done in the following
proposition.

\begin{prop} \label{coerenza1}
$$
D({\cal L})\subset \{u \in {\cal V} \cap W^{2,2}_{loc}:
(1+|x|^\alpha) \Delta u \in L^2_\mu \}
$$
and ${\cal L} u=(1+|x|^\alpha)u$ for $u \in D({\cal L})$. If $\lambda
>0$ and $f \in L^p \cap L^2_\mu$, then
$$
(\lambda-{\cal L})^{-1}f=(\lambda-L)^{-1}f.
$$
\end{prop}
{\sc Proof. } The first part of the proposition easily follows from
local elliptic regularity, testing with any $v \in C_c^\infty$ in
(\ref{forma}). To show the coherence of the resolvents we consider
$f \in C_c^\infty$, $supp\ f \subset B(R)$  and
$u=(\lambda-L_{min})^{-1}f$. Then $u \in D(L)$ solves
$$
\Delta u =\frac{\lambda u}{1+|x|^\alpha}
$$
outside $B(R)$ and is a $C^2$-function. Theorem \ref{risolvente}
implies that $u \in D_{p,max}(L)$ for every $p>N/(N-2)$. If $N>4$,
then $u \in D_{2,max}(L)$ hence $\nabla u \in L^2$, see Proposition
\ref{w2p}, and clearly $u \in L^2_\mu$. This yields $u \in H$ but
not yet $u \in {\cal V}$. To show that $u$ can be approximated with
a sequence of $C_c^\infty$-functions, in the norm of $H$, we fix a
smooth $C^\infty$ function $\eta$ such that $\eta \equiv 1$ in
$B(1)$ and $\eta \equiv 0$ outside $B(2)$ and set $\eta_n(x)=\eta
(x/n)$. Clearly $\eta_n u \to u$ in $L^2_\mu$. Concerning the
gradients we have $\nabla (\eta_n u)=\eta_n \nabla u+u \nabla
\eta_n$. The term $\eta_n \nabla u$ converges to $\nabla u$ in
$L^2$, since $\nabla u \in L^2$ and we have to show that $u \nabla
\eta_n \to 0$ in $L^2$. Since $u \in L^{2^*}$ we can use H\"older's
inequality to deduce
\begin{equation} \label{eq1}
\int_{\R^N}|u|^2 |\nabla \eta_n|^2\, dx \le \frac{C}{n^2}\int_{n \le
|x| \le 2n}|u|^2\, dx \le C_1 \left(\int_{n \le |x|\le 2n}
|u|^{2^*}\right)^{1-\frac{2}{N}}
\end{equation}
 which tends to zero as $n \to
\infty$. This shows $u$ can be approximated with a sequence of
$W^{1,2}$ compactly supported functions and to produce a sequence of
smooth approximants it is now sufficient to use convolutions. Then
$u \in {\cal V}$ and, by integration by parts,
$$
a(u,v)=-(\lambda u-f,v)_{L^2_\mu},
$$
that is $u \in D({\cal L})$ and $\lambda u -{\cal L}u=f$. By
density, this shows the coeherence of the resolvents of $L_{min}$
and ${\cal L}$ for $\lambda >0$, hence of ${\cal L}$ and
$(L,D_{p,max}(L))$, see Theorem \ref{risolvente}. The cases $N=3,4$
require some variants, since $p=2$ does not safisfy the inequality
$p>N/(N-2)$. To show that $u$ belongs to $L^2_\mu$ we use the fact
that $u \in L^p$ with $p=2N/(N-2)$ and therefore
$$
\int_{\R^N}|u|^2\, d\mu \le \left
(\int_{\R^N}|u|^{\frac{2N}{N-2}}\,dx \right )^{1-\frac{2}{N}}\left
(\int_{\R^N}\frac{1}{(1+|x|^\alpha)^{\frac{N}{2}}}\,dx \right
)^{\frac{2}{N}}
$$
Next we show that $\nabla u$ belongs to $L^2$. Since $u \in C^2$ and
$$
\int_{\R^N} \nabla u\nabla v\, dx=\int_{\R^N}(\lambda u-f)v\, d\mu
$$ for every $v \in C_c^\infty$, the same equality holds for every
$v \in W^{1,2}$ having compact support. Taking $v_n=\eta_n u$ we get
$$
\int_{\R^N} \eta_n |\nabla u|^2 dx=\int_{\R^N}(\lambda u-f)\eta_n
u\, d\mu-\int_{\R^N}u\nabla u \cdot \nabla \eta_n\, dx.
$$
Since
$$
\int_{\R^N}u\nabla u \cdot \nabla \eta_n\,
dx=-\frac{1}{2}\int_{\R^N}|u|^2 \Delta \eta_n\, dx
$$ we can proceed
as in (\ref{eq1}) to show that this term tends to zero and hence
$\nabla u \in L^2$. From now one, the proof proceeds as in the case
$N>4$. \qed

We can now strengthen Corollary \ref{indipendenza}.

\begin{prop} \label{indipendenza1}
If $N/(N-2)<p \le \infty$, $2<\alpha<\infty$, then the spectrum of
$(L,D_{p,max}(L))$ lies in $]-\infty, 0[$ and consists of a sequence
$\lambda_n$ of eigenvalues, which are simple poles of the resolvent
and tend to $-\infty$. Each eigenspace is finite dimensional and
independent of $p$.
\end{prop}
{\sc Proof. } Since the resolvents of $(L,D_{p,max}(L))$ and $({\cal
L}, (D({\cal L}))$ are coherent and compact in $L^p$, $L^2_\mu$,
respectively all the assertions except the density of the
eigenfuctions follow from \cite[Proposition 2.2]{arendt} (see also
\cite[Proposition 5.2]{MP} for more details). \qed

Observe that $0$ is in the resolvent set of $L$, since it is
injective. This is clear in $L^p$ or $C_0$ because $\Delta u \in
L^p$ implies $u=0$. However, constant functions are in $H$ if
$\alpha
>N$ and this explains why we work with ${\cal V}$ (constant
functions are never in ${\cal V}$ since ${\cal V}$ embeds into
$L^{2^*}$).

Next we show some methods to estimate the first eigenvalue
$\lambda_1$.

\begin{prop}
The following estimates hold
\begin{equation} \label{Stima1}
\lambda_1\leq
-{\left(\frac{\alpha-2}{2}\right)}^{\frac{2}{\alpha}}\frac{\alpha}{\alpha-2}\frac{(N-2)^2}{4}
\end{equation}
and
\begin{equation} \label{stima2}
\lambda_1\leq -(N-2)\frac{\alpha
\sin\left(\frac{2}{\alpha}\pi\right)}{\pi}.
\end{equation}
\end{prop}

{\sc Proof.}
By Corollary \ref{stima-p}, we obtain
$$\|L^{-1}\|\leq
{\left(\frac{2}{\alpha-2}\right)}^{\frac{2}{\alpha}}\frac{\alpha-2}{\alpha}\frac{p^2}{N(Np-N-2p)}$$
By classical spectral theory then
$$|\lambda_1|\geq {\left(\frac{\alpha-2}{2}\right)}^{\frac{2}{\alpha}}\frac{\alpha}{\alpha-2}\frac{N(Np-N-2p)}{p^2}.$$
The function appearing on the right hand side attaints its maximum for $p=\ds\frac{2N}{N-2}$
 where it reaches the value $${\left(\frac{\alpha-2}{2}\right)}^{\frac{2}{\alpha}}\frac{\alpha}{\alpha-2}\frac{(N-2)^2}{4}.$$
Since  the spectrum of $L$ is independent of $p$ we obtain
(\ref{Stima1}). (\ref{stima2}) is obtained in a similar way from
Proposition \ref{norm-inf}. \qed

Observe that the coefficient
$$
{\left(\frac{\alpha-2}{2}\right)}^{\frac{2}{\alpha}}\frac{\alpha}{\alpha-2}
$$
is always greater than or equal to $1$, and it is $1$ for $\alpha=2,
\infty$. Then (\ref{Stima1}) improves the estimate $\lambda_1 \le
-(N-2)^2/4$ which can be obtained using the classical Hardy
inequality. On the other hand (\ref{stima2}) is better than
(\ref{Stima1}) for large $\alpha$ and small $N$, but worse for
$\alpha$ close to 2 or large $N$.

Since ${\cal L}$ is self-adjoint in $L^2_\mu$, its eigenvalues can
be computed through the Raleigh quotients and, in particular,
$$
-\lambda_1=\min\left \{\int_{\R^N} |\nabla u|^2\, dx :
\int_{\R^N}|u|^2\, d\mu=1 \right \}.
$$
Since
\begin{align*}
 \int_{\R^N}|u|^2\, d\mu &\le \left
(\int_{\R^N}|u|^{2^*}\, dx\right )^{1-2\frac{2}{N}}\left
(\int_{\R^N}\frac{1}{(1+|x|^\alpha)^{N/2}}\, dx \right)^{2/N} \\&
\le C_2^2 \int_{\R^N}|\nabla u|^2\, dx \left
(\int_{\R^N}\frac{1}{(1+|x|^\alpha)^{N/2}}\, dx \right)^{2/N},
\end{align*}
it follows that $ -\lambda_1 \ge \left (C_2^2 L(\alpha)\right)^{-1}
$ where $C_2$ is given by \cite{talenti} and
$$
L(\alpha)=\left (\int_{\R^N}\frac{1}{(1+|x|^\alpha)^{N/2}}\, dx
\right)^{2/N}.
$$
Obesrve also that, when $\alpha \to \infty$, then (formally)
$\lambda_1$ tends to the first eigenvalue of the Dirichlet Laplacian
in the unit ball.

\section*{Appendix}

Here we prove a Hardy-type inequality used throughout the paper.

\begin{prop} \label{hardy}
Let $u\in W^{1,p}(\R^N)$ with compact support, $1< p<\infty$, $\gamma\geq 0$. Then
$$\int_{\R^N}|u|^p|x|^\gamma\leq \left(\frac{p}{\gamma+N}\right)^2\int_{\R^N}|u|^{p-4}|Re(\ov{u}\nabla u)|^2|x|^{\gamma+2}.$$
\end{prop}
{\sc Proof.} Let first $u \in C_c^\infty (\R^N)$ and set $g(t)=u(tx)$. Then
\begin{eqnarray*}
|u(x)|^p &=&|g(1)|^p=-p\int_1^\infty |g|^{p-2}Re\left(\ov{g}\frac{\partial g}{\partial t}\right)dt\\&=&
-p\int_1^\infty|u(tx)|^{p-2}Re(\ov{u}\nabla u(tx))x dt.
\end{eqnarray*}
It follows that
\begin{align*}
\int_{\R^N}&|u(x)|^p|x|^{\gamma}dx \leq p\int_1^\infty dt\int_{\R^N}|u(tx)|^{p-2}|Re(\ov{u}(tx)\nabla u(tx))||x|^{\gamma+1}dx\\&=p\int_1^\infty dt\int_{S_{N-1}}d\omega\int_0^\infty |u(tr\omega)|^{p-2}|Re(\ov{u}(tr\omega)\nabla u(tr\omega)|r^{\gamma+N}dr\\&=
p\int_1^\infty \frac{1}{t^{\gamma+N+1}}dt\int_{S_{N-1}}d\omega\int_0^\infty |u(s\omega)|^{p-2}|Re(\ov{u}(s\omega)\nabla u(s\omega)|s^{\gamma+N}ds\\&
=\frac{p}{\gamma+N}\int_{\R^N}|u(x)|^{p-2}|Re(\ov{u}(x)\nabla u(x)||x|^{\gamma+1}dx.
\end{align*}
By density this inequality holds for every $u \in W^{1,p}(\R^N)$ having compact support.
At this point H\"{o}lder's inequality yields
$$
\int_{\R^N}|u(x)|^p|x|^{\gamma}dx \leq \frac{p}{\gamma+N}\left(\int_{\R^N}|u(x)|^p|x|^{\gamma}dx\right)^\frac{1}{2}\left(\int_{\R^N}|u(x)|^{p-4}|Re(\ov{u}(x)\nabla u(x)||x|^{\gamma+2}dx\right)^\frac{1}{2}.
$$
\qed

\end{document}